\documentstyle[draft, amscd, syntonly, psfig, 
amssymb, 12pt]{amsart}
\chardef\bslash=`\\
\newcommand{\ntt}{\series m\shape n\tt}
\newcommand{\cs}[1]{{\protect\ntt\bslash'1}}
\newcommand{\opt}[1]{{\protect\ntt'1}}
\newcommand{\env}[1]{{\protect\ntt'1}}
\makeatletter
\def\verbatim{\interlinepenalty\@M \@verbatim
\leftskip\@totalleftmargin\advance\leftskip2pc
\frenchspacing\@vobeyspaces \@xverbatim}
\makeatother
\newtheorem{theo}{Theorem}[section]

\newtheorem{lem}[theo]{Lemma}
\newtheorem{defi}[theo]{Definition}
\newtheorem{remark}[theo]{Remark}
\newtheorem{corol}[theo]{Corollary}

\def \proof {{\bf Proof$\colon$}\ }
\def \Z{{\bf Z}}

\def\nt{$n$-trivializer}

\begin{document}

\title[]{Knot adjacency and satellites}
\author{E. Kalfagianni}
\author{ X.-S. Lin}
\address[]{Department of Mathematics , Michigan State
University,
E. Lansing, MI, 48823}
\email[]{kalfagia@@math.msu.edu}
\address[]{Department of Mathematics , 
University of California,
Riverside, CA,  92521}
\email[]{xl@@math.ucr.edu}
\thanks{The research of the first author is partially supported
by the NSF through grant DMS-0104000. The research of the second
author is partially supported
by the Overseas Youth Cooperation
Research Fund of NSFC and by the NSF through grant DMS-0102231.}
\begin{abstract} A knot $K$ is called $n$-adjacent to the unknot, if $K$ 
admits a projection containing $n$ {\em generalized crossings}
such that changing any $0<m\leq n$ of them yields a
projection of the unknot. We show that a non-trivial satellite knot
$K$ is $n$-adjacent to the unknot, for some $n>0$, if and only
if it
is $n$-adjacent to the unknot
in any companion solid torus.
In particular,
every {\em model knot} of $K$ is
$n$-adjacent to the unknot. Along the way of proving these results, we
also show that
2-bridge knots of the form $K_{{p/q}}$,
where ${\displaystyle {p/q}}=[2q_1,2q_2]$ for some $q_1,q_2\in \Z$,
are precisely those knots that have
genus one and are $2$-adjacent to the unknot.
\vskip 0.12in



\end{abstract}

\maketitle


\section{Introduction}

The development of the theory of finite type knot invariants 
has led to the notion of {\em $n$-triviality} which is a multiplex unknotting operation.
This notion was introduced independently by Gussarov
and Ohyama (\cite{kn:g}, \cite{kn:O}). Roughly speaking,
a knot is $n$-trivial if it can be unknotted in
$2^n-1$ different ways by  {\it multiple
crossing changes}. 
The research in this paper is motivated by
the following question: If a non-trivial satellite knot $K$ is $n$-trivial
is there 
a companion  torus of $K$ that is disjoint from all the crossing changes
that exhibit $K$ as $n$-trivial?
In this paper, we are concerned with a stronger version of $n$-triviality where each set 
of multiple crossing changes is taken to be a set of {\em twist crossings}
on two strings of the knot. A knot with this stronger $n$-triviality is called
$n$-{\it adjacent} to the unknot. 
For knots which are $n$-adjacent to the unknot, using results of 
Lackenby(\cite{kn:la}) and Scharlemann-Thompson 
(\cite{kn:st1}, \cite{kn:st2}), we
obtain an affirmative answer to the aforementioned question.
In fact, we show that the generalized crossings involved can be taken to be
disjoint from any companion torus of $K$. As a consequence, we obtain that
if a non-trivial satellite knot
$K$ is $n$-adjacent to the unknot then it is
$n$-adjacent to the unknot
in any companion solid torus.
In particular,
any {\em model knot} of $K$ is
$n$-adjacent to the unknot. Along the way of proving these results, we
also characterize
2-bridge knots of the form $K_{{p/q}}$,
where ${\displaystyle {p/q}}=[2q_1,2q_2]$ for some $q_1,q_2\in \Z$,
precisely as those knots that have
genus one and are $2$-adjacent to the unknot.
\smallskip

A generalized crossing of order $q\in \Z$
on an embedding of a knot $K$ is a set $C$
of $|q|$ twist crossings
on two strings that inherit
opposite orientations from any orientation of $K$.
If $K'$ is obtained from $K$ by changing
all the crossings in $C$ simultaneously, we will say that
$K'$ is obtained from $K$ by a generalized crossing change (see Figure 1). 
In particular, if $|q|=1$, $K$ and $K'$ differ by an ordinary crossing change
while if $q=0$ we have $K=K'$. Note that a generalized crossing change can be achieved 
by ${\displaystyle {1\over q}}$-surgery on a crossing circle, which is
an unknotted curve that bounds an embedded disc $D\subset S^3$
such
that $K$ intersects ${\rm int}(D)$ exactly twice with
zero algebraic intersection number.
{\vspace{.03in}}
\begin{figure}[ht]
\centerline{\psfig{figure=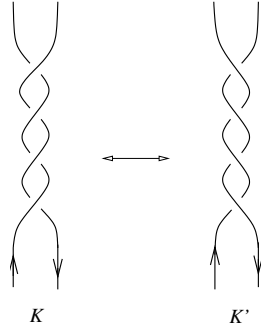,height=1.7in, clip=}}
\caption{The knots $K$ and $K'$ differ by a generalized crossing change
of order $q=-4$.}
\end{figure}

\begin{defi} \label{defi:adj}
We
will say that $K$ is $n$-adjacent to the unknot, for some
$n\in {\bf N}$,
if $K$ admits an embedding containing $n$ generalized crossings
such that changing any $0<m\leq n$ of them yields an
embedding of  the unknot. A collection of crossing circles
corresponding to these crossings is called an \nt.
If all the generalized crossings used
have order
$+1$ or $-1$ (i.e. they are ordinary crossings), we will say that
$K$ is simply $n$-adjacent to the unknot. An \nt\  that shows $K$ to be
simply $n$-adjacent to the unknot will be called a simple \nt.
\end{defi}

\begin{remark}{\rm Let $V$ be a solid torus in $S^3$ and suppose that a
knot $K$ is embedded in $V$. Throughout the paper, we will use the term
``$K$ is $n$-adjacent to the unknot in $V$" to mean the following:
There exists an embedding of $K$ in $V$
that contains $n$ generalized crossings
such that changing any $0<m\leq n$ of them unknots $K$ in $V$.} 
\end{remark}

To state our result
recall that if $K$ is
a non-trivial satellite with
{\em companion knot} ${\hat K}$  and {\em model knot} $P$ then:
i) ${\hat K}$ is non-trivial;
ii) $P$ is geometrically essential in a standardly embedded
solid torus $V_1\subset S^3$;  and iii) there is a homeomorphism
$h: V_1\longrightarrow  V:=h(V_1),$
such that $h(P)=K$ and $\hat K$ is the core of $V$.

\begin{theo} \label{theo:main} Let $K$ be a non-trivial satellite knot
and let $V$ be any companion solid torus of $K$.
Then,
$K$ is
$n$-adjacent to the unknot, for some $n>0$, if and only if
it is $n$-adjacent to the unknot in $V$.
\end{theo}

As a consequence of Theorem \ref{theo:main} we have the following:

\begin{corol} \label{corol:model} A non-trivial satellite 
knot $K$ is
$n$-adjacent to the unknot, for some $n>0$, if and only if
any model knot of $K$ is 
$n$-adjacent to the unknot in the standard solid torus $V_1$.
\end{corol}

The paper is organized as follows: In Section 2 we summarize some results
that are used in the proof of Theorem \ref{theo:main}. 
In Section 3 we study satellite knots
of winding number zero that are $n$-adjacent to the unknot.
In Section 4 we prove Theorem \ref{theo:main} and work out some corollaries.
In Section 5, we show that a knot $K$ of genus one is 2-adjacent to 
the unknot iff it is
a 2-bridge knot of the form $[2q_1,2q_2]$ for some $q_1,q_2\in\Z$.

Note that a weaker version of Theorem \ref{theo:main} is 
generalized to a broader class of $n$-trivial knots in
\cite{kn:k}.
\smallskip

\noindent{\bf Acknowledgments:} We would like to thank the referee
and Ying-Qing Wu for their careful reading and critical remarks 
of the previous version of this paper.

\section{Preliminaries} 
In this section we summarize some results that
will be used in the proof of Theorem \ref{theo:main}. We begin with the 
following
theorem that summarizes some results from \cite{kn:st1} and \cite{kn:la}.
Part a) of the theorem is stated as Corollary 3.2 in \cite{kn:st1}.
Part b) is stated as Corollary 4.4 of \cite{kn:st1} and also
follows from 
Theorem 1.4(b) and Proposition 2.5 of
\cite{kn:la}.

\begin{theo} \label{theo:stl}Let $K, K'$ be knots
that differ by a generalized crossing change of order $q\in \Z$. Let $L$ be a crossing circle for $K$
corresponding to this generalized crossing and that $K$ is a non-trivial satellite.

a)\ If  $|q|>1$ and  ${\rm genus}(K')\leq {\rm genus}(K)-1$ then
any companion torus of $K$ can be isotoped in $S^3\setminus K$ to be disjoint from $L$.

b)\ If $|q|=1$ and ${\rm genus}(K')\leq {\rm genus}(K)-2$ then
any companion torus of $K$ can be isotoped in $S^3\setminus K$ to be disjoint from $L$.
\end{theo}

To continue, note that if $L$ is an \nt \ for a knot $K$,
since the linking number of $K$ with each component of $L$ is zero,
$K$ bounds a Seifert surface in the complement of $L$. We will need the following lemma:

\begin{lem} \label{lem:taut} Let $L$ be an \nt\ of a
knot $K$. 
Suppose that $S$ is a Seifert surface bounded by $K$ in the complement
of $L$ and such that among all such surfaces
$S$ has minimum genus. Then,
${\rm genus}(S)={\rm genus}(K)$. 
\end{lem}

\proof For simple \nt s the lemma is stated as 
Theorem 4.1 in \cite{kn:hl}; the proof relies on a result of Gabai
(Corollary 2.4 of \cite{kn:ga}). The argument for
general \nt's is essentially the same. The details are given 
in the proof of Theorem 3.1 of \cite{kn:kl}. 
\qed
\smallskip

We will apply Theorem \ref{theo:stl} to the case when $K'$ is the unknot.
In this case, when ${\rm genus}(K)\geq 2$, for every $q$, we can isotopy 
a companion torus of $K$ to be disjoint from the crossing
circle $L$. When ${\rm genus}(K)=1$ and $|q|=1$, Theorem \ref{theo:stl} can
not be used anymore. Nevertheless, the following lemma guarantees that the
same conclusion still holds in this situation. 

\begin{lem} \label{lem:plum}Let $K, K'$ be knots
that differ by a generalized crossing of order $q\in \Z$. Let $L$ be the 
corresponding crossing circle. Suppose that
${\rm genus}(K')<{\rm genus}(K)$ and let $K_0$ denote the 2-component
link obtained by smoothing $C$ in a way consistent with the orientation
of $K$. Then, there exist  Seifert surfaces $\Sigma$ and $\Sigma_0$ of
maximal Euler characteristic 
for $K$ and $K_0$ respectively, $\Sigma\cap L=\emptyset$, 
such that $\Sigma$ is obtained from
$\Sigma_0$ by plumbing on an unknotted annulus with a 
$(2, 2q)$-torus link as its boundary and $L$ as one of its 
small linking circles.
\end{lem}

Here by a {\it small linking circle} of an annulus $A$ embedded in $S^3$,
we mean an unknot which bounds a disk $D$ such that $D\cap A$ is a
proper non-boundary parallel arc in $A$. 
\smallskip

\proof For $|q|=1$ the lemma is stated as 
Proposition 3.1 in \cite{kn:st2}. The proof of that proposition uses
Theorem 1.4 (of \cite{kn:st2}) that gives a relation between the 
Euler characteristics
of the triple $(K,K', K_0)$. Theorem 6.4.3 of \cite{kn:kaiser}
states that the same relation holds when $|q|>1$.
Using this, the arguments used in the proof of Proposition 3.1 
in \cite{kn:st2}
go through to give the lemma in the case that $|q|>1$ (see also
the proof of Theorem 6.4.2 of \cite{kn:kaiser}). \qed
\smallskip

\begin{corol}\label{corol:all} Let $K$ be a non-trivial satellite knot 
which can be unknotted by a single generalized
crossing change, and $L$ be the corresponding crossing circle. Then
any companion torus of $K$ can be isotoped in $S^3\setminus K$ to be disjoint 
from $L$.  
\end{corol}
\proof The case when ${\rm genus}(K)\geq 2$ is covered by Theorem 
\ref{theo:stl}.
So we assume ${\rm genus}(K)=1$. Let $\Sigma$ be the genus one Seifert 
surface of 
$K$ claimed to exist in Lemma \ref{lem:plum}. Then $\Sigma$ is the 
plumbing of two annuli $A_1$ and $A_2$. 
One of them, say $A_1$, is unknotted with $L$ as its small linking circle. 
Thus,
$K$ is contained in a torus $T=\partial N$, where $N$ is a tubular 
neighborhood of
$A_2$. Since $A_1$ is unknotted, we may assume that $A_1\subset N$. 
It is not hard to see that $T$ is the innermost companion torus of $K$:
every other companion torus of $K$ can be isotoped to contain the solid 
torus $N$
in one side.
Thus, we can isotope every companion torus of $K$ to be disjoint 
from $L$. \qed  

\section{Satellite knots with zero winding number}

Throughout this section, we suppose that $K$ is a non-trivial satellite 
knot and $V$ is a companion solid torus of $K$, such that the winding number 
of $K$ in $V$ is zero.
\smallskip

\subsection{A technical lemma}
In this subsection we prove a technical lemma which will play
a key role in our discussion in the next subsection.
\smallskip

Let $S$ be a minimal genus Seifert surface of $K$. We assume that the 
intersection of
$S$ and $T=\partial V$ is transverse and the number of components 
of $S\cap T$ is minimal. Denote by $M_1$ and $M_2$ the closures of components of $(S^3\setminus K)\setminus T$ in $S^3\setminus K$, respectively, with $M_2$
a compact 3-manifold.
Let $\alpha$ be a proper arc on $S$.

\begin{lem}\label{lem:red} Suppose there is an isotopy of $S^3$, fixing 
$K$ pointwise, which brings $\alpha$ to an arc $\alpha'$ in $V$, then we 
can isotopy $\alpha$ on the
surface $S$, relative to $\partial S=K$, to a proper arc $\alpha''$ in $V$.
\end{lem}

\proof We may assume that

(1) $S\cap T$ is a collection of disjoint parallel copies of an 
essential simple closed curve on $T$;

(2) every component of $S\cap M_i$ is incompressible and boundary incompressible in $M_i$. 
\smallskip

Since these points follow from well known facts in 3-dimensional manifold
topology, we only give a brief explanation.
Point (1) follows immediately from the fact that $T$ is incompressible
in the complement of $K$. To see (2), first by the incompressibility of $S$ and $T$ in $S^3\setminus K$, 
it is easy to deduce that each component of $S\cap M_i$ is incompressible in $M_i$. If there is an essential boundary compressing disk $D$ for a component of $S\cap M_i$ in $M_i$, then $D\cap T$ must be an arc whose end points lie on different components of $S\cap T$. Thus we may isotopy $S$ to reduce
$|S\cap T|$, which would contradict the assumption that $|S\cap T|$ is 
minimal.   
\smallskip

Now let $\alpha$ and $\alpha'$ be as in the lemma.  Up to isotopy on $S$, relative to $\partial S=K$, we
may assume that $\alpha$ intersects each component of $S\cap T$ in
essential arcs. We will to show that this assumption will force $\alpha$ to be disjoint from $T$. Let $f: D^2\rightarrow S^3$ be a path homotopy from $\alpha$ to
$\alpha'$ with $f(\text{Int}(D))$ disjoint from $K$.  Since $T$ is
incompressible, we may assume that $f^{-1}(T)$ is a set of proper arcs on
$D^2$.  Note that all endpoints of $f^{-1}(T)$ are on $\alpha$ because
$\alpha'$ is disjoint from $T$.  Thus we can choose a component
$\beta$ of $f^{-1}(T)$ which is outermost in the sense that it cuts
off a subdisk $D_1$ in $D^2$ whose interior is disjoint from
$f^{-1}(T)$, and $\gamma = f(\partial D_1\setminus\text{Int}(\beta))$ is a subarc of $\alpha$. Since the interior of $D_1$ is disjoint from
$f^{-1}(T)$, $\gamma$ is an proper arc on a component $A$ of $S\cap M_i$, 
which is
essential by the above assumption. However, this contradicts the following lemma and the fact that $A$ is
boundary incompressible in $M_i$. So we conclude that $\alpha$ can be isotoped on $S$, relative to $\partial S=K$, to be disjoint from $T$. This proves the lemma. \qed

\begin{lem} Let $F$ be an incompressible and boundary incompressible surface
in a 3-manifold $M$ with $\partial M$ incompressible.  Then there is no non-closed proper essential curve $a$ on $F$ that is homotopic to a curve $b$ on
$\partial M$ relative to $\partial a$.
\end{lem}

\proof Consider the double of $F$ in the double of $M$, denoted by $\hat
F$ and $\hat M$, respectively.  By an innermost-circle outermost-arc
argument one can easily show that $\hat F$ is incompressible in $\hat
M$.  On the other hand, the double of a homotopy from $a$ to $b$ would give
rise to a null homotopy disk for the double of $a$.  Since the double of
$a$ is an essential curve on $F$, this contradicts the fact that an
incompressible surface is $\pi_1$-injective. \qed

\subsection{Finding an \nt\ in a companion solid torus} 
We can now have the the following
lemma, which will allows us to find
an \nt\ for $K$ in any companion solid torus.

\begin{lem} \label{lem:torus} Let $K$ be  a non-trivial satellite
and let $V$ be a companion solid torus of $K$.
Suppose that the winding number of $K$ in $V$
is zero. If $K$ is $n$-adjacent to the unknot, for some $n>0$,
then, there exists an \nt\  for $K$ that lies in
$V$.
\end{lem}

\proof
Let $L:=\cup_{i=1}^n L_i$
be an $n$-trivializer of $K$
and let $D_1,\ldots, D_n$ be crossing discs
bounded by $L_1,\ldots,  L_n$, respectively. 
Let $S$ be a Seifert surface
for $K$ in the complement
of $L$ that has minimum genus. Then by Lemma \ref{lem:taut},
$S$ is also a minimal genus Seifert surface of $K$. 
We may isotope $S$ so that each $S\cap {\rm int}(D_i)$ is the union of an arc
$\alpha_i$ and several closed components. The arc $\alpha_i$ 
is properly embedded on $S$.
Since $S$ is incompressible in the 
complement of $L$, after an isotopy we can arrange so that  
$S\cap D_i$ contains no closed curves that are 
inessential on $D_i$. 
Thus each
closed component of $S\cap D_i$ has to be parallel to $L_i$ on $D_i$.
By replacing $L_i$ with the closed component of $S\cap D_i$
that is innermost on $D_i$, we may assume that $S\cap D_i=\alpha_i$. 
Since twisting along $\alpha_i$ unknots
$K$, it follows that $\alpha_i$ must be essential on $S$.
Furthermore, the arcs $\alpha_1,\dots,\alpha_n$ are disjoint from each other.
\smallskip

By Corollary \ref{corol:all}, for each $L_i$, we can isotope the 
torus $T=\partial V$ in the complement of $K$ to $T'$ such that
$T'\cap L_i=\emptyset$. Assume that $T'$ intersects the disk
$D_i$ transversely. Since $T'$ is disjoint from $L_i=\partial D_i$,
each component of $T'\cap D_i$ is a simple closed curve in $D_i$. If
a component of $T'\cap D_i$ bounds a disk in $D_i$ which contains
only one point in $K\cap D_i$, we would have the winding number of
$K$ in $V$ to be $\pm1$. So every component of $T'\cap D_i$ either
bounds a disk in $D_i$ which is disjoint from $K\cap D_i$ or
bounds a disk in $D_i$ which contains $K\cap D_i$. In either cases,
a further isotopy of $T'$ in the complement of $K$ will remove this component
of $T'\cap D_i$. The reversed isotopy in the complement of $K$
from this $T'$ to $T$ then will bring the arc $\alpha_i$ into $V$.
Thus, by Lemma \ref{lem:red}, we can isotopy each $\alpha_i$ on the 
minimal genus Seifert surface 
$S$, relative to $\partial S=K$, to a proper arc $\alpha_i'$ in $V$. 
\smallskip

On $S$, let $\alpha$ and $\beta$ be two proper 1-submanifolds whose 
intersection is transverse.
Suppose that there is an isotopy
of $S$ that reduces the geometric intersection
$|\alpha\cap\beta|$. Then there will be 
a disk $D$ on $S$ such that $D\cap(\alpha\cup\beta)=\partial D$, 
$D\cap\alpha$ and $D\cap\beta$
are subarcs in the interior of $\alpha$ and $\beta$, respectively.
This is a well-known fact (see, for example, Proposition 3.10 in
\cite{kn:thurston} or Lemma 3.1 in \cite{kn:hass}).
\smallskip

We apply
this fact to $\{\alpha_1,\dots,\alpha_n\}$ and 
$S\cap T=C_1\sqcup\cdots\sqcup C_r$ (see the proof of Lemma \ref{lem:red}). 
Since $\alpha_1$ can be made disjoint from $S\cap T$ by an isotopy of 
$S$ relative to 
$\partial S$, we find a disk $D$ between $\alpha_1$ and $S\cap T$ as
described above. 
We then use this $D$ to define an isotopy of $S$ relative to 
$\partial S$ to remove 
a pair of intersection points of $\alpha_1$ and $S\cap T$. 
This isotopy will not 
increase the intersection points of the other $\alpha_i$'s with $S\cap T$. 
And it will also keep $\alpha_i$'s disjoint. So inductively,
we have an isotopy on $S$ relative to $\partial S$, which brings the entire 
disjoint collection of proper arcs $\{\alpha,\dots,\alpha_n\}$ to a 
disjoint collection
of proper arcs $\{\alpha'_1,\dots,\alpha'_n\}$ in $V$.
\smallskip

Finally, we construct a small disk $D'_i$ in $V$ whose intersection 
with $S$ is 
$\alpha'_i$, for each $i$, and they are disjoint from each other. 
Let 
$L'_i=\partial D'_i\subset V$. Since $L'=\cup_{i=1}^n L'_i$ is isotopic 
to $L$ in the 
complement of $K$, it is an \nt\ for $K$ that lies in $V$. \qed

\section{The proof of the main result} Here we finish the proof of Theorem \ref{theo:main} and Corollary \ref{corol:model}.
\smallskip

{\bf Proof of Theorem \ref{theo:main}:} The ``if" direction of the statement 
is clear. To prove
the ``only if" direction, suppose that  $K$ is a non-trivial satellite
that is $n$-adjacent to the unknot. Let $V$ be any companion solid torus 
of $K$. Let 
$\hat K$ denote 
the core of $V$ and set $T:=\partial V$. 
We need the following:
\smallskip

{\it Claim.} The winding number of $K$ in $V$ is zero.
\smallskip

{\it Proof of Claim.}  
By Corollary \ref{corol:all},
there exists a component $L_1\subset L$ that can be isotoped
to be disjoint from $T$.
Let $D_1$ be a crossing disc
bounded by $L_1$.  After an isotopy in the complement of $K$, $D_1\cap T$ 
will consist
of a collection of curves, none of which
bounds a disc in $D_1$ in the complement of $K$.
Let $C$ be a component of $D_1\cap T$. If $C$ is boundary parallel on $D_1$
then it can be eliminated by an isotopy in the complement of $K$ so 
that $L_1$ 
is still disjoint from $T$. If all components of $D_1\cap T$ are 
boundary parallel,
we will have $D_1$ disjoint from $T$ after an isotopy in the complement 
of $K$.
Then $D_1$ is contained in $V$. Since a 
satellite with non-zero winding number cannot be unknotted
by crossing changes in ${\rm Int}(V)$, we conclude every
component of $D_1\cap T$
bounds a disc on $D_1$ that contains
exactly one point of $D_1 \cap K$.  
Since $K$ was assumed to be
a non-trivial satellite we conclude that $K$ is a
composite knot and $T$ is the follow-swallow torus.
But then the crossing change realized by $L_1$ occurs within
a summand of $K$ and it cannot unknot $K$. This contradicts
the fact that $L_1$ is part of an $n$-trivializer and it finishes the proof
of the claim.
\smallskip

Let us now finish the proof of the theorem. The claim above allows us to assume that
the winding number of $K$ in $V$ is zero.
By Lemma \ref{lem:torus},
$K$ admits an \nt \ $L'$ in $V$.
Now each of the surgeries
along the sublinks of $L'$
that unknot $K$ must turn it into a knot that is isotopically trivial in 
Int$(V)$.
For, otherwise the knot obtained from K after any of these surgeries will 
still have $\hat K$ as a companion and it can't be the unknot. Thus $K$ is $n$-adjacent 
to the unknot in $V$. \qed
\smallskip

{\bf Proof of Corollary \ref{corol:model}:} Let $P$ be any model of 
$K$ in a standard solid torus $V_1\subset S^3$ and let $h: \longrightarrow S^3$ 
the satellite embedding.
If
$P$ is $n$-adjacent to the unknot in $V_1$
and 
$L\subset V_1$ is  an \nt\ then $h(L)$ is an \nt \ 
for $K$ in $V$. Conversely, by Theorem \ref{theo:main} and its proof,
if $K$ is  $n$-adjacent to the unknot then any \nt, say $L$,  can be 
isotoped into
$V:=h(V_1)$ as an \nt of $K$ in $V$.
But then the crossing circles $h^{-1}(L_1), \ldots, h^{-1}(L_n)$ 
form an \nt \ for
$P$ in $V_1$. \qed
\smallskip

There exist many criteria in terms of the finite type knot invariants or
polynomial invariants that detect $n$-adjacency to the unknot.
For example,
in
\cite{kn:ak} it is shown that
if a knot is $n$-adjacent to the unknot, for some 
$n\geq 3$, then all the finite type invariants of order $<2n-1$
and the Alexander
polynomial are trivial. More recently,
criteria that detect simple 2-adjacency to the unknot
were obtained by N. Askitas and A. Stoimenow (\cite{kn:as}) in terms of the HOMFLY
polynomial, and
by the second named author of this paper and  Z. Tao in terms of the Kauffman
polynomial.
Due to the computational complexity of the invariants involved,
these criteria become 
harder to test for knots that are non-trivial satellites.
The results of this paper,
reduce the
problem of deciding whether 
a non-trivial satellite $K$
is $n$-adjacent to the unknot
to deciding the same problem for a
model knot of $K$. In particular, we have the following:

\begin{corol} Suppose that $P$ is a knot that is not $n$-adjacent
to the unknot. Then, no satellite that is modeled on $P$ is
$n$-adjacent to the unknot.
\end{corol}

\proof This follows from Theorem \ref{theo:main} and the fact that if 
$P$ is $n$-adjacent to the unknot in the solid torus then it is 
$n$-adjacent to the unknot in $S^3$. \qed

\section{Knots of genus one}

We finish this paper by taking a look at genus one knots that are
$n$-adjacent to the unknot, for $n>1$. In fact, we will obtain
a characterization of knots of genus one which are 2-adjacent to
the unknot.
\smallskip

Consider 2-bridge knots of the form $K_{p/q}$, where $p/q=[2q_1,2q_2]$
in Conway's notation (see, for example, \cite{kn:bz}). Such a knot is
formed by plumbing two unknotted, $2q_1$ and $2q_2$ twisted annuli, and
taking the boundary of the resulting genus one surface. Obviously,
this is a genus one knot which is 2-adjacent to the unknot. The orders
of the two generalized crossing changes are $q_1$ and $q_2$, respectively.
\smallskip

\begin{theo}\label{theo:genusone} A genus one knot $K$ is 2-adjacent
to the unknot if and only if $K=K_{p/q}$, $p/q=[2q_1,2q_2]$, for some
integers $q_1,q_2$.
\end{theo}

It is clear that we only need to prove the \lq\lq only if" part. So we
suppose $K$ is a genus one knot and it is 2-adjacent to the unknot.
Let $L=L_1\cup L_2$ be a 2-trivializer of $K$ of order $(q_1,q_2)$.
By Lemma \ref{lem:taut}, we have a Seifert surface $S$ of $K$ in the
complement of $L$ and the genus of $S$ is one. We may assume that the
crossing disks $D_1,D_2$, with $\partial D_1=L_1$ and $\partial D_2=L_2$,
intersect $S$ along essential proper arcs $\alpha_1,\alpha_2$,
respectively.

\begin{lem}\label{lem:parallel} {\rm (See also \cite{kn:hl}.)}
The arcs $\alpha_1,\alpha_2$ are not parallel on $S$.
\end{lem}
\proof If $\alpha_1$ and $\alpha_2$ were parallel to each other on $S$,
$L_1$ and $L_2$ would cobound an annulus in the complement of $K$. Then
perform both $1/q_1$-surgery on $L_1$ and $1/q_2$ surgery on $L_2$ would
be the same as doing $1/(q_1+q_2)$-surgery on $L_1$ or $L_2$. Since $K$
is nontrivial, we would have two distinct surgeries on $L_1$ under which
$S$ does not remain of minimal genus.

Since a twist along $L_1$ unknots $K$, it follows that the 3-manifold
$M:=S^3\setminus \eta(K\cup L_1)$ is irreducible.
$S$ gives rise to a properly embedded surface in $M$ that minimizes the
Thurston norm in its homology class. Corollary 2.4 of \cite{kn:ga},
applied to $M$ and  $T:=\partial \eta(L_1)$,
implies that there can be at most one
Dehn filling of $T$ (or equivalently at most one surgery along $L_1$)
under which $S$ doesn't remain a minimum genus surface for $K$.
This contradiction finishes the proof of the lemma.
\qed
\smallskip

Next we use Lemma \ref{lem:plum} to $K$ and $L_1$. This lemma gives us
another genus one Seifert surface $\Sigma$ of $K$ in the form of the
plumbing of annuli $A_1$ and $A_2$, such that $A_1$ is unknotted,
$2q_1$-twisted, and its small linking circle is $L_1$. Without loss of
generality, we may assume that $D_1\cap\Sigma=\alpha_1$. So we may cut
open both $S$ and $\Sigma$ along $\alpha_1$. For $\Sigma$, we get the
annulus $A_2$ from this surgery. For $S$, we get another annulus $A'_2$
from this surgery. The annuli $A_2$ and $A'_2$ have the same boundary
$K_0$, and we may assume that they are disjoint. Thus $T=A_2\cup A'_2$
is a torus which bounds a solid torus in $S^3$. Assume that the core circle
of $A_2$ (and the core circle of $A'_2$) is a $(m,l)$ curve on $T$, where
$l\geq 0$ is the winding number in the longitude direction and $m$ is 
the winding number in the meridian direction on $T$. 
\smallskip

The arc $\alpha_2$ on $A'_2$ should have its end points on different boundary 
components of $A'_2$ by Lemma \ref{lem:parallel}. We may pick a possible
$\alpha_2$ and all other possible $\alpha_2$'s are obtained by Dehn twist 
along $A'_2$. See Figure 2.

{\vspace{.03in}}
\begin{figure}[ht]
\centerline{\psfig{figure=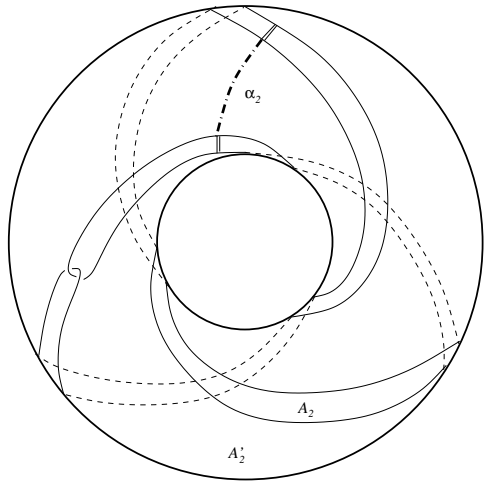,height=2.5in, clip=}}
\caption{A possible position for the arc $\alpha_2$ on the annulus $A'_2$.}
\end{figure}

Let us perform a generalized crossing change at $L_1$, which will unhook the 
clasp seen on $A_2$ in Figure 2. We then shrink the two separated clasp 
ends along $A_2$ until they meet the ends of $\alpha_2$ (the places on
$A_2$ marked by double lines in Figure 2). Denote the subarc of the
core circle of $A_2$ between the double line marks, which does not run 
through the clasp, by $\beta$. 
Then we may get a simple closed curve $J=\alpha_2\cup\beta$. The curve $J$ 
intersects the core circle of $A'_2$ only once. Let $J$ be a $(a,b)$ curve
on $T$. We have $|mb-la|=1$. 
\smallskip

If we perform generalized crossing changes at both $L_1$ and $L_2$, $K$ will 
be changed to the unknot. This is possible only when $(a,b)=(0,\pm1)$ 
or $(a,b)=(\pm1,0)$. Otherwise, $J$ would have a non-zero framing in the 
solid torus $V$ bounded by $T$ and generalized crossing changes at both 
$L_1$ and $L_2$ would change $K$ into a $[2r,2s]$ knot in $V$ for
$rs\neq0$. Such a knot can not be unknotted in $S^3$.  
\smallskip

When $(a,b)=(0,\pm1)$, $V$ has to be unknotted in $S^3$. Notice that we 
must have $m=\pm1$. If $l=0$, the knot $K$ would be trivial.
So we may assume that $l\geq1$.
Thus, we can have one possible choice of $\alpha_2$ as shown in 
Figure 3.  
Any other choices of $\alpha_2$ are obtained by applying a power of the 
Dehn twist of along the core of $A'_2$ to this particular $\alpha_2$. From 
this fact, we see that this particular $\alpha_2$ is the only one with
$(a,b)=(0,\pm1)$. In Figure 3, we can see that the arc $\alpha_2$ can be 
isotoped to the arc $\gamma$ in $A'_2$. Furthermore, $\gamma$ and 
an arc on $A_2$, which is the intersection of $A_2$ with its small 
linking disk, cobound a disk whose interior is disjoint from $T$. Thus 
the generalized crossing change at $L_2$ is the same as a generalized 
crossing change at a small linking circle of $A_2$. This implies that
$K$ is a 2-bridge knot of the form $[2q_1,2q_2]$. 
\smallskip

When $(a,b)=(\pm1,0)$, we have $l=1$ and $m\neq0$. 
We can argue as in the previous case: First, we find a unique choice of 
$\alpha_2$. This choice of $\alpha_2$ will 
force $V$ to be unknotted in $S^3$. And then the knot $K$ will be
a 2-bridge knot of the form $[2q_1,2q_2]$. This finishes the proof of 
Theorem \ref{theo:genusone}. \qed

{\vspace{.03in}}
\begin{figure}[ht]
\centerline{\psfig{figure=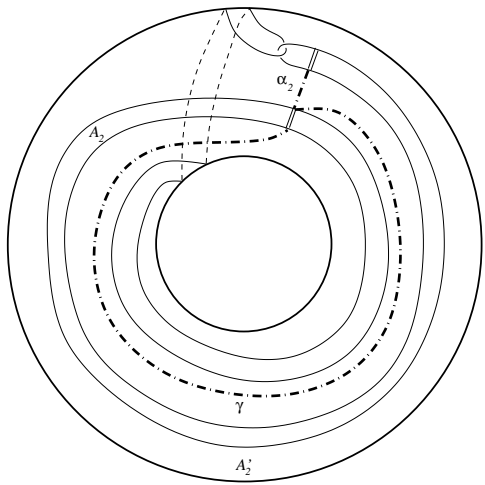,height=2.5in, clip=}}
\caption{This is a knot of the form $[2q_1,2q_2]$.}
\end{figure}

\begin{corol} The only genus one
knots
that are simply 2-adjacent to the unknot are the two trefoils and the figure
eight.
\end{corol}

\proof This corollary corresponds to the case of $q_1=\pm1$ and $q_2=\pm1$
of Theorem \ref{theo:genusone}.\qed

Note that, as observed by T. Stanford, if a knot $K$
is simply 2-adjacent to the unknot, then we have $a_2(K)=0$ or $\pm1$,
where
$a_2$ is the second coefficient of the Alexander-Conway polynomial.
Using this observation, we see that the knot $5_2$ is not simply
2-adjacent
to the unknot since $a_2(5_2)=2$.
By Theorem \ref{theo:genusone}, $5_2$ is not 2-adjacent to the unknot since it is
the 2-bridge knot $[2,3]$.
Apparently, no method is known to detect this using
knot invariants.

\end{document}